\documentclass[reqno,12pt]{amsproc}
  \usepackage{latexsym}
  \usepackage[all]{xy}
  \usepackage{amsfonts}
  \usepackage{amsthm}
  \usepackage{amsmath}
  \usepackage{amssymb}
  \usepackage{pifont}

  \def\sw#1{{\sb{(#1)}}}

  \def\suc#1{{\sp{(#1)}}}

  \def\proof{{\sl Proof.~~}}

  \def\endproof{\hbox{$\sqcup$}\llap{\hbox{$\sqcap$}}\medskip}
  \def\<{{\langle}}
  \def\>{{\rangle}}

  \def\eps{\varepsilon}

  \def\note#1{{}}

  \def\note#1{}

  \def\cD{{\mathcal D}}

  \def\beq{\begin{equation}}
  \def\eeq{\end{equation}}

  \def\id{\mathrm{id}}

  \def\ot{{\otimes}}

     \def\1{\mathbf{1}}

\def\act{\!\cdot\!}
\def\lqmod#1{{}_{#1}\mathcal{Q}\mathcal{M}}
\def\rqmod#1{\mathcal{Q}\mathcal{M}_{#1}}
\def\Lqmod#1{{}^{#1}\mathcal{Q}\mathcal{M}}
\def\Rqmod#1{\mathcal{Q}\mathcal{M}^{#1}}

\def\rcosmash{\!\blacktriangleright\!\!<\!}
\def\lsmash{\!>\!\!\!\triangleleft}

\def\k{\Bbbk}

  \newcounter{zlist}
  \newenvironment{zlist}{\begin{list}{(\arabic{zlist})}{
  \usecounter{zlist}\leftmargin2.5em\labelwidth2em\labelsep0.5em
  \topsep0.6ex
  \parsep0.3ex plus0.2ex minus0.1ex}}{\end{list}}

  \newcounter{blist}

  \newcounter{rlist}

   \newcounter{alist}

  \renewcommand{\subjclassname}{\textup{2010} Mathematics Subject
        Classification}

\def\stac#1{\raise-.2cm\hbox{$\stackrel{\displaystyle\otimes}{\scriptscriptstyle{#1}}$}}

\def\cten#1{\raise-.2cm\hbox{$\stackrel{\displaystyle\widehat{\otimes}}
{\scriptscriptstyle{#1}}$}}

  \def\Label#1{\label{#1}\ifmmode\llap{[#1] }\else
  \marginpar{\smash{\hbox{\tiny [#1]}}}\fi}
  \def\Label{\label}

  \newtheorem{proposition}{Proposition}[section]
  \newtheorem{lemma}[proposition]{Lemma}
  \newtheorem{corollary}[proposition]{Corollary}
  \newtheorem{theorem}[proposition]{Theorem}

  \theoremstyle{definition}
  \newtheorem{definition}[proposition]{Definition}

  \newtheorem{example}[proposition]{Example}

  \theoremstyle{remark}
  \newtheorem{remark}[proposition]{Remark}

  \theoremstyle{definition}
  \swapnumbers

\begin{document}

 \title{Actions of Hopf quasigroups}
 \author{Tomasz Brzezi\'nski}
 \address{{\bf Tomasz Brzezi\'nski: }Department of Mathematics, Swansea University,
  Singleton Park,
  Swansea SA2 8PP, U.K.}
  \email{T.Brzezinski@swansea.ac.uk}
  
 \author{Zhengming Jiao}
 \address{{\bf Zhengming Jiao: } Department of Mathematics,
Henan Normal University,
Xinxiang 453002, Henan,
Peoples Republic of China } 
 \email{zmjiao@henannu.edu.cn}

\date{May 2010}

  \begin{abstract}
Definitions of actions of Hopf quasigroups are discussed in the context of Long dimodules and smash products. In particular, Long dimodules are defined for Hopf quasigroups and coquasigroups, and solutions to Militaru's $\cD$-equation are constructed. A necessary compatibility condition between action and multiplication of a Hopf quasigroup acting on its quasimodule Hopf quasigroup for a smash product construction is derived.\\

\noindent{\bf \subjclassname :} {16T05; 16T15}
  \end{abstract}
  \maketitle

\section{Introduction}
{\em Hopf quasigroups} and {\em Hopf coquasigroups} were  introduced recently in \cite{KliMaj:Hop} in order to capture the quasigroup features of the (algebraic) 7-sphere. These are generalisations of Hopf algebras that are not required to be (co)associative. The aim of this paper is to discuss possible definitions of actions of Hopf quasigroups and coactions of Hopf coquasigroups. Since Hopf quasigroups are not required to be associative, the natural compatibility between the multiplication and action should be different from the associative law. It seems, however, that there is no universal compatibility condition that could be imposed. What relationship between the multiplication and action should be required depends on the application in mind. For example, in \cite{Brz:qua}, where Hopf modules for Hopf quasigroups were defined, a particular compatibility condition was requested that stems naturally from the relationship between quasigroups (or loops) and Hopf quasigroups. In this paper we use the same condition for the definition of Long dimodules and we prove that a different (stronger) condition is needed for construction of smash products of Hopf quasigroups. 

A {\em Long dimodule} over a Hopf algebra $H$ is a vector space $M$ with a left $H$-action and a right $H$-coaction that is required to be left $H$-linear (with respect to the left $H$-module structure on $H\otimes M$ induced by the multiplication in $H$). Long dimodules were introduced in \cite{Lon:Bra} in the context of Brauer groups, and were used in \cite{Mil:Lon} to solve a class of nonlinear equations termed  {\em $\mathcal{D}$-equations} which we describe presently. Given a vector space $M$, a linear endomorphism $R: M\ot M\to M\ot M$ is said to be a solution of the $\mathcal{D}$-equation if
\begin{equation}\label{eq.d}
R^{12}\circ R^{23} = R^{23}\circ R^{12},
\end{equation}
where $R^{12} = R\ot \id, \; R^{23} = \id\ot R :  M\ot M\ot M\to M\ot M\ot M$. 
In the first part of this paper we introduce Long-dimodules for Hopf quasigroups and Hopf coquasigroups and show that they are a source of solutions to the $\mathcal{D}$-equation.

Smash or cross products of Hopf quasigroups and coproducts of Hopf coquasigroups were introduced in \cite{KliMaj:Hop} as the `quasi' versions of their Hopf algebra predecessors \cite{Mol:sem}.  The construction of a smash product involves  two Hopf quasigroups, one acting on the other. The action is assumed to be associative. Since multiplication in a Hopf quasigroup is not associative this assumption might appear artificial. In the second part of this paper we show that, in fact, this assumption is {\em almost} necessary (it is truly necessary in the case of a bijective antipode). More precisely, we study the smash product construction without assuming that the action is associative, and then derive the necessary condition  for the existence of the antipode, Eq.~\eqref{modass}, which turns out to be the associativity of the action up to the antipodal operation. A similar analysis is then carried out for smash coproducts of Hopf coquasigroups.

All algebras and coalgebras are over a field $\k$. Unadorned tensor product symbol represents the tensor product of $\k$-vector spaces. 

\section{Preliminaries on Hopf (co)quasigroups}\label{sec.prelim}
\setcounter{equation}{0}
The aim of this section is to recall the definition of a Hopf quasigroup and a Hopf coquasigroup from \cite{KliMaj:Hop}. 
Let $H$ be a vector space that is a unital (not necessarily associative) algebra with product $\mu: H\ot H\to H$ and unit $1: \k \to H$, and a counital (not necessarily coassociative) coalgebra with coproduct $\Delta: H\to H\ot H$ and counit $\eps: H\to \k$ that are algebra homomorphisms.

$H$ is called a {\em Hopf quasigroup} provided $\Delta$ is coassociative and there exists a linear map $S: H\to H$ such that
\begin{equation}\label{quasi1}
\mu\circ (\id \ot \mu) \circ (S\ot \id\ot \id )\circ (\Delta \ot \id ) = \eps \ot \id  = \mu\circ (\id  \ot \mu) \circ (\id \ot S\ot \id )\circ (\Delta \ot \id )
\end{equation}
and
\begin{equation} \label{quasi2}
\mu\circ (\mu \ot \id ) \circ (\id \ot \id \ot S)\circ (\id \ot \Delta) = \id \ot \eps = \mu\circ (\mu \ot \id ) \circ (\id \ot S\ot \id )\circ (\id \ot \Delta).
\end{equation}

$H$ is called a {\em Hopf coquasigroup} provided $\mu$ is associative and there exists a linear map $S: H\to H$ such that
\begin{equation}\label{coq1}
(\mu \ot \id ) \circ (S\ot \id \ot \id )\circ (\id \ot \Delta)\circ \Delta = 1\ot \id  = (\mu \ot \id ) \circ (\id \ot S\ot \id )\circ (\id \ot \Delta)\circ \Delta
\end{equation}
and
\begin{equation}\label{coq2}
(\id  \ot \mu) \circ (\id \ot \id \ot S)\circ (\Delta \ot \id )\circ \Delta = \id \ot 1 = (\id  \ot \mu) \circ (\id \ot S\ot \id )\circ (\Delta \ot \id )\circ \Delta .
\end{equation}

We use Sweedler notation for coproduct: for all $h\in H$, $\Delta (h) = h\sw 1\ot h\sw 2$ (summation implicit).  Thus, in terms of the Sweedler notation, the Hopf quasigroup conditions \eqref{quasi1}--\eqref{quasi2} come out as,
$$
S(h\sw 1)(h\sw 2g) = h\sw 1(S(h\sw 2) g) = (gh\sw 1)S(h\sw 2) = (gS(h\sw 1))h\sw 2 = g\eps(h),
$$
for all $g,h\in H$. Dually, the Hopf coquasigroup conditions  \eqref{coq1}--\eqref{coq2} come out as 
$$
S(h\sw 1)h\sw 2\sw 1\ot h\sw 2\sw 2 = h\sw 1S(h\sw 2\sw 1)\ot h\sw 2\sw 2  = 1\ot h
$$
and
$$
h\sw 1\sw 1 \ot h\sw 1\sw 2 S(h\sw 2) = h\sw 1\sw 1 \ot S(h\sw 1\sw 2) h\sw 2 = h\ot 1.
$$
Note that since $\Delta$ is not coassociative in this case, the standard Sweedler's relabeling rules cannot  be used. 

As for standard Hopf algebras, the map $S$  is called an {\em antipode}. 
It is proven in \cite{KliMaj:Hop} that the antipode is antimultiplicative and anticomultiplicative and it immediately follows from (any of) equations \eqref{quasi1}--\eqref{coq2} that, for all $h\in H$,
$
S(h\sw 1)h\sw 2 = h\sw 1 S(h\sw 2) = \eps(h)1,
$
i.e.\ $S$ enjoys the standard antipode property.

\section{Quasimodules and Long dimodules for Hopf quasigroups}\label{sec.Long} \setcounter{equation}{0}
The aim of this section is to introduce Long dimodules over Hopf quasigroups and Hopf coquasigroups and to prove that they provide solutions to the Militaru $\cD$-equation \eqref{eq.d}.

\subsection{Quasimodules over Hopf quasigroups}\label{sec.qmod}

First recall from \cite{Brz:qua}  the definition of a quasimodule over a Hopf quasigroup and a quasicomodule over a Hopf coquasigroup.

\begin{definition}[\cite{Brz:qua}]\label{def.quasim} 
Let $H$ be a Hopf quasigroup,  a {\em left $H$-quasimodule} is a pair $(M, \rho_M)$, where $M$ is a $\k$-vector space, $\rho_M: H\ot M\rightarrow M$, $h\ot m \mapsto h\act m$, is a linear map such that, for all $m\in M$ and $h\in H$,
\begin{equation}\label{qm2}
1\act m=m, \qquad 
h_{(1)}\act (S(h_{(2)})\act m)=\varepsilon (h)m=S(h_{(1)})\act (h_{(2)}\act m).
\end{equation}
\end{definition}

A morphism of left $H$-quasimodules is defined as a left $H$-linear map, i.e.\ a $\k$-linear map $f:M\to N$ such that, for all $m\in M$, $h\in H$, $f(h\cdot m) = h\cdot f(m)$. The category of left $H$-quasimodules over a Hopf quasigroup $H$ is denoted by $\lqmod H$. 

\begin{remark}\label{rem.loop}
The definition of a quasimodule is intended to reflect that a Hopf quasigroup is a linearisation of a {\em loop}. Recall from \cite{Alb:qua} that a {\em loop} is a set $L$ with a binary operation $L\times L\to L$, $(a,b)\mapsto ab$ that has a neutral element $e\in L$ (i.e.\ $ae=ea =a$) and such that, for all $a\in L$, there exists $a^{-1}\in L$ satisfying, for all $b\in L$,
$$
a(a^{-1}b) = a^{-1}(ab) = (ba^{-1})a = (ba)a^{-1} = b.
$$
If $L$ is a loop then the $\k$-algebra $H=\k L$ spanned by the elements of $L$ is a Hopf quasigroup with unit $e$ and multiplication induced from that in $L$, the comultiplication $\Delta(a) =a\ot a$, counit $\eps(a) =1$ and antipode $S(a) = a^{-1}$, for all $a\in L$. 

When discussing actions of loops on sets it is natural to require that a left $L$-set is a set $X$ with an operation $L\times X\to X$, $(a,x)\mapsto a\cdot x$, that is compatible with the unit and inverses in $L$, i.e.\
$$
e\act x = x,\qquad a^{-1}\act (a\act x) = a\act (a^{-1}\act x) = x.
$$
Then $X$ is an $L$-set if and only if $\k X$ is a $\k L$-quasimodule.
\end{remark}
\begin{remark} \label{rem.ex}
Let $H$ be a Hopf quasigroup.

(1)  $H$ is a left $H$-quasimodule with the multiplication as action and so is the tensor product $H\ot M$, where $M$ is an $H$-quasimodule. The $H$-action on $H\ot M$ is given by the multiplication in $H$ and by identity on $M$. Note, however,  that the action $\rho_M:H\ot M\to M$ is not a morphism of quasimodules.  Note, further,  that an associative left $H$-module $M$ can be viewed as a left $H$-quasimodule. In this case the action is a morphism in $\lqmod H$.

(2) The field $\k$ is a right $H$-quasimodule with trivial action $h\act \alpha= \varepsilon (h)\alpha$. Furthermore, if  $M,N\in \lqmod H$,  then $M\ot N$ is also an object in $\lqmod H$ with the {\em diagonal action} 
$$
h\act (m\ot n)=h_{(1)}\act m\ot h_{(2)}\act n.
$$
This makes the category of quasimodules into a monoidal category $(\lqmod H, \ot, \k)$.
\end{remark}

The category $\rqmod H$ of right $H$-quasimodules is defined similarly. Dually,

\begin{definition}[\cite{Brz:qua}]\label{def.coquasim} 
Let $H$ be a Hopf coquasigroup,  a {\em right $H$-quasicomodule} is a pair $(M, \rho^M)$, where $M$ is a $\k$-vector space and  $\rho^M: M\rightarrow M\ot H$, $m\mapsto m^{(0)}\ot m^{(1)}$ (summation implicit), is a linear map such that, for all $m\in M$,
\begin{equation}\label{cqm2}
m^{(0)}\varepsilon (m^{(1)})=m, \quad 
m^{(0)(0)}\ot m^{(0)(1)}S(m^{(1)})=m\ot 1=m^{(0)(0)}\ot S(m^{(0)(1)})m^{(1)}.
\end{equation}
\end{definition}

The category of right $H$-quasicomodules over a Hopf coquasigroup $H$ with $H$-colinear maps as morphisms is denoted by $\Rqmod H$. The category $\Lqmod H$ of left $H$-quasicomodules is defined similarly. 

By considerations dual to those in Remarks~\ref{rem.loop}--\ref{rem.ex}, examples of right $H$-co\-quasi\-co\-modules include $(H,\Delta)$ and $(\k, 1 )$. Also, if $(M,\rho^M)$ is a quasicomodule, $\rho^M$ is not a morphism in $\Rqmod H$ from $(M,\rho^M)$ to $(M\ot H, \id\ot \Delta)$. Furthermore,  $(\Rqmod H, \ot, \k)$ is a monoidal category, where the coaction on $M\ot N$ is given by
$$
\rho ^{M\ot N}(m\ot n)=m^{(0)}\ot n^{(0)}\ot m^{(1)}n^{(1)}.
$$

\subsection{Long dimodules over Hopf quasigroups}

\begin{definition}\label{def.ldim} 
Let $H$ be a Hopf quasigroup. A {\em Long $H$-dimodule} is a triple $(M, \rho_M, \rho^M)$, where $(M, \rho_M)$ is a left $H$-quasimodule, $(M, \rho^M)$ is a unital and coassociative right $H$-comodule such that the following   compatibility condition holds 
\begin{equation}\label{ldm}
\rho^M\circ \rho_M=(\rho_M\otimes H)\circ (H\otimes \rho^M).
\end{equation}
\end{definition}

Writing  $h\act m$ for $\rho_M(h\otimes m)$ and using the Sweedler notation for coaction $\rho^M(m)=m^{(0)}\otimes m^{(1)}$, the conditions \eqref{ldm} read explicitly, for all $h\in H$ and $m\in M$,
\begin{equation}\label{ldm1}
(h\act m)^{(0)}\otimes (h\act m)^{(1)}=h\act m^{(0)}\otimes m^{(1)}.
\end{equation}

The category of Long $H$-dimodules over a Hopf quasigroup $H$ with $H$-linear $H$-colinear maps
is denoted by $_{\underline H}{\mathcal L}^H$.

\begin{example}
Let $L$ be a loop  and consider a family of $L$-sets $\{X_a \; |\; a\in L\}$. Define $M:= \bigoplus _{a\in L} \k X_a$. Then $M$ is a Long dimodule over a Hopf quasigroup $\k L$, provided each of the $X_a$ is equipped with the $\k L$-coaction $x\mapsto x\ot a$ and  $\bigoplus _{a\in L} \k X_a$ is equipped with the resulting direct sum action and coaction; compare \cite[Example~3.2(1)]{Mil:Lon}.
\end{example}

\begin{lemma} \label{def.ldmp} 
Let $H$ be a Hopf quasigroup, $(M,\rho_M,\rho^M) \in  {}_{\underline H}{\mathcal L}^H$. Then, for all $m\in M$ and $h\in H$,
\begin{equation}\label{ldmp1}
\rho^M(S(m^{(1)})\act m^{(0)})= S({m^{(1)}}_{(2)})\act m^{(0)}\otimes {m^{(1)}}_{(1)}
\end{equation}
and
\begin{equation}\label{ldmp2}
\rho^M(h\act (S(m^{(1)})\act m^{(0)}))=h\act (S({m^{(1)}}_{(2)})\act m^{(0)})\otimes {m^{(1)}}_{(1)}.
\end{equation}
\end{lemma}

\proof Equation \eqref{ldmp1} follows by the compatibility condition \eqref{ldm1} and by the coassociativity of the coaction $\rho^M$. Then equation \eqref{ldmp2} is a simple consequence of \eqref{ldm1} and \eqref{ldmp1}. \endproof

The following propositions are the Hopf-quasigroup versions of \cite[Example~3.2 \& Remark~3.3]{Mil:Lon}. 
\begin{proposition} \label{def.ldmqmt}  
Let $H$ be a Hopf quasigroup. Then, for any left $H$-quasimodule $(M, \rho _M)$,  $M\otimes H$ is a Long  $H$-dimodule with action $\rho _{M\ot H} = \rho_M\ot \id$ and coaction $\rho ^{M\ot H} = \id\ot\Delta$.
The functor $\bullet \otimes H: \lqmod H \rightarrow _{\underline H}\!\!{\mathcal L}^H$ is the right adjoint of the forgetful functor $ _{\underline H}{\mathcal L}^H\rightarrow \lqmod H$.
\end{proposition}

\proof By the unitality of action $\rho_M$,  $1\act (m\ot k)=m\ot k$. For all $m\in M, k, h\in H$, 
$$
h_{(1)}\act (S(h_{(2)})\act (m\ot k))=h_{(1)}\act (S(h_{(2)})\act m)\ot k=\varepsilon (h)m\ot k,
$$
where the first equality follows by the definition of $\rho _{M\ot H}$ and the second and third ones by  equations \eqref{qm2}.  The coassociativity of comultiplication $\Delta$ and the counitality of $\varepsilon $ imply that $(M\ot H, \rho ^{M\ot H})$ is a right $H$-comodule. Finally, combining the definitions of $\rho _{M\ot H}$ and $\rho ^{M\ot H}$ we obtain
$$
\rho^{M\ot H} (h\act (m\ot k))=\rho^{M\ot H} (h\act m\ot k)
=h\act m\ot k_{(1)}\ot k_{(2)}
=h\act (m\ot k)^{(0)}\ot (m\ot k)^{(1)}.
$$
Therefore, the compatibility condition \eqref{ldm1} between $\rho _{M\ot H}$ and $\rho ^{M\ot H}$ holds, and $M\ot H$ is a Long $H$-dimodule. 

The unit $\eta$ and the counit $\sigma$ of adjunction are defined as follows. For all Long $H$-dimodules $(M,\rho_M,\rho^M)$, $\eta_M = \rho^M$, and, for all quasimodules $(N,\rho_N)$, $\sigma_N = \id \ot \eps: N\ot H\to N$. Note that $\eta_M$ is a morphism of Long $H$-dimodules by the coassociativity of $\rho^M$ and the compatibility condition \eqref{ldm}. The triangular identities follow by the counitality of $\rho^M$.
\endproof

\begin{proposition}\label{def.ldmcmt}   
Let $H$ be a Hopf quasigroup. For any right $H$-comodule $(M, \rho ^M)$,  $H\otimes M$ is a Long  $H$-dimodule by the action $\rho _{H\ot M} = \mu\ot\id$ and coaction $\rho ^{H\ot M} = \id\ot\rho^M$.
\end{proposition}

\proof 
 By the definition of $\rho _{H\ot M}$ and $\rho ^{H\ot M}$,  the Hopf quasigroup stucture of $H$ and the right $H$-comodule structure of $M$ make $H\ot M$ a left $H$-quasimodule and a right $H$-comodule, respectively. We need only to check the  compatibility  between $\rho _{H\ot M}$ and 
 $\rho ^{H\ot M}$. As a matter of fact, 
$$
\rho^{H\ot M} (h\act (k\ot m))=\rho^{H\ot M} (hk\ot m)
=hk\ot m^{(0)}\ot m^{(1)}
=h\act (k\ot m)^{(0)}\ot (k\ot m)^{(1)}.
$$
This completes the proof. 
\endproof

\begin{remark}
Note the asymmetry  in statements of Propositions~\ref{def.ldmqmt} and \ref{def.ldmcmt}. In the latter we stopped short from claiming that the functor $H\ot \bullet: {\mathcal M}^H\rightarrow {}_{\underline H}{\mathcal L}^H$ and  the forgetful functor $ _{\underline H}{\mathcal L}^H\rightarrow {\mathcal M}^H$ form an adjoint pair. For a Long $H$-dimodule $(M, \rho_M, \rho^M)$ the expected counit of this adjunction $\sigma_M$ would be given by the action $\rho_M$, which, as explained in Remark~\ref{rem.ex}(1), is not a morphism of quasimodules, hence not a morphism in $_{\underline H}{\mathcal L}^H$. Note that the expected unit $\eta_N = \id \ot 1$ is a comodule map, and that $\sigma$ and $\eta$ satisfy the triangular identities by the unitality of $\rho_M$.
\end{remark}
As explained in Remark~\ref{rem.ex}(2) the field $\k$ is a left $H$-quasimodule by the counit $\eps$,  and it is a right $H$-comodule with the  coaction given by the unit of $H$. It follows from Propositions~\ref{def.ldmqmt} and \ref{def.ldmcmt} that $H$ is a Long $H$-dimodule in two different ways: with multiplication in $H$ and the trivial coaction $\id\ot 1 $, and with trivial action $\id\ot \eps$ and comultiplication $\Delta$. 
Generalising this observation we obtain the following corollary of Propositions~\ref{def.ldmqmt} and \ref{def.ldmcmt}.

\begin{corollary} \label{def.qmtoldm}   
(1) Let $H$ be a Hopf quasigroup and let $(M, \rho _M)$ be a left $H$-quasimodule. Then $(M, \rho _M, \rho ^M)$, where $\rho ^M: M\rightarrow M\ot H$, $m\mapsto m\ot 1$,  is a Long  $H$-dimodule.

(2) Let $(M, \rho ^M)$ be a right comodule over  a Hopf quasigroup $H$. Then $(M, \rho _M, \rho ^M)$, where $\rho _M: H\ot M\rightarrow M$, $h\ot m \mapsto \varepsilon (h)m$,  is a Long  $H$-dimodule.
\end{corollary}

Since $(\k,\eps)$ is an $H$-quasimodule, Corollary~\ref{def.qmtoldm} implies that $(\k,\eps, 1 )$ is a Long $H$-dimodule which is the first step in establishing that $(_{\underline H}{\mathcal L}^H, \ot, \k)$ is a monoidal category; see \cite[Remark~3.3(2)]{Mil:Lon} for the Hopf algebra case. The second step, the definition of a tensor product of Long $H$-dimodules is contained in the following

\begin{proposition}\label{def.ldmts}   
Let $H$ be a Hopf quasigroup, $(M, \rho _M, \rho^M)$ and $(N, \rho_N, \rho^N)$ be Long $H$-dimodules. Then $M\ot N$ is  a Long $H$-dimodule with the diagonal action and coaction, for all $m \in M$, $n\in N$, $h\in H$, 
$$
\rho_{M\ot N}(h\ot m\ot n):=h_{(1)}\act m\ot h_{(2)}\act n \quad 
\mbox{and}  \quad
\rho^{M\ot N}(m\ot n):=m^{(0)}\ot n^{(0)}\ot m^{(1)}n^{(1)}.
$$
\end{proposition}

\proof 
Since $1\in H$ is a grouplike element and actions $\rho_M$ and $\rho_N$ are unital, also $\rho_{M\ot N}$ is unital. To check the quasimodule property, take any
 $m, \in M$ and $n\in N$, and compute
\begin{align*}
h_{(1)}\act (S(h_{(2)})\act (m\ot n))
&= h_{(1)}\act (S(h_{(3)})_{(1)}\act m)\ot h_{(2)}\act (S(h_{(3)})_{(2)}\act n)\cr
&= h_{(1)}\act (S(h_{(4)})\act m)\ot h_{(2)}\act (S(h_{(3)})\act n)
=\varepsilon (h)m\ot n,
\end{align*}
by the anticomultiplicativity of the antipode and the quasimodule properties of $M$ and $N$.
Similarly,  $S(h_{(1)})\act (h_{(2)}\act (m\ot n))=\varepsilon (h)m\ot n$. Therefore, $(M\ot N, \rho_{M\ot N})$ is a left $H$-quasimodule.  $(M\ot N, \rho ^{M\ot N})$ is a right $H$-comodule by the coassociativity of coactions and multiplicativity of the coproduct in $H$. It remains to check that the compatibility condition \eqref{ldm1} holds.  For all $m\in M$, $n\in N$ and $h\in H$,
\begin{align*}
\rho^{M\ot N}(h\act (m\ot n))&=\rho^{M\ot N}(h_{(1)}\act m\ot h_{(2)}\act n)\\
&=(h_{(1)}\act m)^{(0)}\ot (h_{(2)}\act n)^{(0)}\ot (h_{(1)}\act m)^{(1)}(h_{(2)}\act n)^{(1)}\\
&=h_{(1)}\act m^{(0)}\ot h_{(2)}\act n^{(0)}\ot m^{(1)}n^{(1)}
=h\act (m\ot n)^{(0)}\ot (m\ot n)^{(1)},
\end{align*}
where the Long $H$-dimodule compatibility conditions for $M$ and $N$ have been used to obtain the second equality.
\endproof

Finally, we construct a solution to the Militaru $\cD$-equation \eqref{eq.d}.

\begin{proposition}\label{def.D-eq}    
Let $H$ be a Hopf quasigroup and let  $(M, \rho _M, \rho^M)$ be a Long $H$-dimodule. Then the  map 
$$
R_{(M, \rho _M, \rho^M)}: M\ot M\rightarrow M\ot M,
\;\;\;
m\ot n\mapsto  n^{(1)}\act m\ot n^{(0)}, 
$$
is a solution of the ${\mathcal D}$-equation.
\end{proposition}

\proof 
Write $R$ for $R_{(M, \rho _M, \rho^M)}$. Then, by the Long $H$-dimodule compatibility condition \eqref{ldm1}, for all $l, m, n\in M$,
\begin{align*}
R^{12}\circ R^{23}(l\ot m\ot n)&=R^{12}(l\ot n^{(1)}\act m\ot n^{(0)})
=(n^{(1)}\act m)^{(1)}\act l\ot (n^{(1)}\act m)^{(0)}\ot n^{(0)}\\
&=m^{(1)}\act l\ot n^{(1)}\act m^{(0)}\ot n^{(0)}
=R^{23}(m^{(1)}\act l\ot m^{(0)}\ot n) \\
&=R^{23}\circ R^{12}(l\ot m\ot n), 
\end{align*}
 i.e.\ $
R^{12}\circ R^{23}=R^{23}\circ R^{12},
$ as required. \endproof

\subsection{Long dimodules over Hopf coquasigroups}

Considerations regarding Long dimodules over Hopf coquasigroups are dual to these regarding Hopf quasigroups, so we merely outline main observations without giving detailed proofs.

\begin{definition}\label{def.ldimc} 
Let $H$ be a Hopf coquasigroup. A Long $H$-dimodule is a triple $(M, \rho_M, \rho^M)$, where $(M, \rho_M)$ is a unital and associative left $H$-module and $(M, \rho^M)$ is a right $H$-quasicomodule such that the compatibility condition \eqref{ldm} holds.
\end{definition}

Note that we require a Long $H$-dimodule over a Hopf coquasigroup to satisfy the same compatibility condition as in the case of a Hopf quasigroup as  \eqref{ldm} is self-dual, i.e.\ it is invariant under simultaneous interchanging actions with coactions, reversing the order of composition and flipping tensor products (this last operation ensures that $H$ acts on the left and coacts on the right of $M$).
The category of right Long $H$-dimodules over a Hopf coquasigroup $H$ with $H$-linear $H$-colinear maps as morphisms will be denoted by $_H{\mathcal L}^{\overline H}$.  Dually to Lemma~\ref{def.ldmp} we have

\begin{lemma} \label{def.ldmcp} 
Let $H$ be a Hopf coquasigroup, $M\in  {}_H{\mathcal L}^{\overline H}$. Then, for all $m\in M$,
\begin{equation}\label{ldmcp1}
\rho^M(S(m_{(1)})\act m_{(0)})= S({m^{(1)}}_{(2)})\act m_{(0)}\otimes {m^{(1)}}_{(1)}
\end{equation}
and, for all $m\in M, h\in H$,
\begin{equation}\label{ldmcp2}
\rho^M(h\act (S(m_{(1)})\act m_{(0)}))=h\act (S({m^{(1)}}_{(2)})\act m_{(0)})\otimes {m^{(1)}}_{(1)}
\end{equation}
\end{lemma}

As explained at the end of Section~\ref{sec.qmod}, the coaction of a quasicomodule is not a morphism of quasicomodules. On the other hand, in view of the compatibility condition \eqref{ldm} the action is $H$-colinear and is linear by the associativity, hence it is a morphism of long $H$-dimodules. Therefore, Proposition~\ref{def.ldmqmt} and Proposition~\ref{def.ldmcmt} dualise to, respectively,

\begin{proposition}   
For any right quasicomodule $(M, \rho ^M)$ over a Hopf coquasigroup $H$,  $H\otimes M$ is a Long  $H$-dimodule with action $\rho _{H\ot M} = \mu\ot\id$ and coaction $\rho ^{H\ot M} = \id\ot\rho^M$.
The functor $H\ot \bullet: \Rqmod H\rightarrow _H\!\!{\mathcal L}^{\overline H}$ is the left adjoint of the forgetful functor $ _H{\mathcal L}^{\overline H}\rightarrow \Rqmod H$.
\end{proposition}
\proof We only indicate the forms of the unit $\eta$ and counit $\sigma$ of this adjunction. For any $H$-comodule $(N, \rho^N)$, $\eta_N = 1 \ot \id  : N\to H\ot N$, and, for all Long-dimodules $(M,\rho_M, \rho^M)$, $\sigma_M = \rho_M: H\ot M \to M$.
\endproof

\begin{proposition}   
For any left module $(M, \rho _M)$ over a Hopf coquasigroup $H$, $M\otimes H$ is a Long  $H$-dimodule with action $\rho _{M\ot H} = \rho_M\ot \id$ and coaction $\rho ^{M\ot H} = \id\ot\Delta$.
\end{proposition}

Given a Hopf coquasigroup $H$, every $H$-module can be made into a Long $H$-dimodule and every $H$-quasicomodule can be made  into a Long $H$-dimodule  as in Proposition~\ref{def.qmtoldm}. In particular $H$ itself is a Long $H$-dimodule in two ways. The category ${}_H{\mathcal L}^{\overline H}$ is a monoidal category with unit $\k$ and tensor product as in Proposition~\ref{def.ldmts}. The formula in Proposition~\ref{def.D-eq} with $M\in {}_H{\mathcal L}^{\overline H}$ gives a solution of Militaru's $\cD$-equation \eqref{eq.d}.

\section{Module algebras for Hopf quasigroups and smash products}\label{sec.smash} \setcounter{equation}{0}
Smash or cross products for Hopf quasigroups and coproducts for Hopf coquasigroups were introduced in \cite{KliMaj:Hop}. The former involve two Hopf quasigroups such that one acts on the other respecting its algebra and coalgebra structures, i.e.\ $H$ and $A$ are Hopf quasigroups such that $A$ is an $H$-module  Hopf quasigroup. The action of $H$ on $A$ is assumed to be associative even though $H$ itself is not an associative algebra. The aim of this section is to re-examine the smash product construction from \cite{KliMaj:Hop} not assuming that $A$ is an $H$-module algebra from the onset, but rather that $A$ is an $H$-{\em quasi}module Hopf quasigroup. We then conclude that the associativity  of the $H$-action up to an application of the antipode is a necessary condition for the smash product construction. In this way we slightly improve and clarify results of \cite{KliMaj:Hop}. In the second part of this section we discuss the smash coproduct of Hopf coquasigroups. 

\subsection{Module algebras  and smash products for Hopf quasigroups}\label{subsec.smash}

We start with the definition of a quasimodule Hopf quasigroup.

\begin{definition}\label{def.qmalg}
Let $H$ be a Hopf quasigroup,  a (not necessarily associative) algebra $A$ is called a left {\em $H$-quasimodule algebra} if $A$ is a left $H$-quasimodule and 
\begin{equation}\label{qmalg}
(h_{(1)}\act a)(h_{(2)}\act b)=h\act(ab),  \qquad h\act 1_A=\varepsilon (h)1_A,
\end{equation}
for all $h, g\in H, a, b, c\in A$. 

A coalgebra $C$ is a left {\em $H$-quasimodule coalgebra} if $C$ is a left $H$-quasimodule and 
\begin{equation}\label{qmcoalg}
\Delta (h\act c)=h_{(1)}\act c_{(1)}\ot h_{(2)}\act c_{(2)}, \qquad \varepsilon (h\act c)=\varepsilon (h)\varepsilon (c),
\end{equation}
for all $h\in H, c\in C$. 

A Hopf quasigroup $A$ is called a {\em left  $H$-quasimodule Hopf quasigroup} if it is both a left $H$-quasimodule algebra and coalgebra (by the same $H$-action). 
\end{definition}

Noting that the proof of \cite[Lemma~4.13]{KliMaj:Hop} does not rely on the associativity of the action we can state

\begin{lemma}\label{def.qmS}
The antipode of a left $H$-quasimodule Hopf quasigroup is $H$-linear.
\end{lemma}

The following theorem, which gives the necessary condition for the construction in \cite[Proposition~4.14]{KliMaj:Hop}, is the main result of this section.

\begin{theorem}\label{def.smash}
Let $H$ be a Hopf quasigroup,  $A$ a left $H$-quasimodule Hopf quasigroup such that, for all $h\in H$ and $a\in A$, 
\begin{equation}\label{cocommu}
h_{(1)}\ot h_{(2)}\act a=h_{(2)}\ot h_{(1)}\act a.
\end{equation}
Then the following statements are equivalent:
\begin{zlist}
\item
There is a smash product Hopf quasigroup $A\lsmash  H$ built on $A\ot H$ with tensor product coproduct, counit and unit, and 
\begin{equation}\label{smashmulti}
(a\ot h)(b\ot g)=a(h_{(1)}\act b)\ot h_{(2)}g,
\end{equation}
\begin{equation}\label{smashanti}
S(a\ot h)=S(h_{(2)})\act S(a)\ot S(h_{(1)}),
\end{equation}
for all $a, b\in A, g, h\in H$.
\item 
For all $g, h\in H$ and $a\in A$, 
\begin{equation}\label{modass}
g\act (S(h)\act a)=(gS(h))\act a
\end{equation}
\end{zlist}
\end{theorem}

\proof $(1)\Rightarrow (2)$ If $A\lsmash  H$ forms a Hopf quasigroup in the given way, then the second of the Hopf quasigroup conditions \eqref{quasi2} reads
$$
((b\ot g)S((a\ot h)_{(1)}))(a\ot h)_{(2)}
=\varepsilon (a)\varepsilon (h)b\ot g.
$$
Developing the left hand side of this condition we obtain
\begin{align*}
((b\ot g)&S((a\ot h)_{(1)}))(a\ot h)_{(2)}
=((b\ot g)S(a_{(1)}\ot h_{(1)}))(a_{(2)}\ot h_{(2)})\\
&\stackrel{\eqref{smashanti}}{=} ((b\ot g)(S(h_{(1)})_{(1)}\act S(a_{(1)})\ot S(h_{(1)})_{(2)}))(a_{(2)}\ot h_{(2)})\\
&\stackrel{\eqref{smashmulti}}{=}(b(g_{(1)}\act (S(h_{(1)})_{(1)}\act S(a_{(1)})))\ot g_{(2)}S(h_{(1)})_{(2)})(a_{(2)}\ot h_{(3)})\\
&\stackrel{\eqref{smashmulti}}{=} (b(g_{(1)}\act (S(h_{(1)})_{(1)}\act S(a_{(1)}))))((g_{(2)}S(h_{(1)})_{(2)})\act a_{(2)})\ot (g_{(3)}S(h_{(1)})_{(3)})h_{(2)}\\
&\;\;{=}\;\; (bS(g_{(1)}\act (S(h_{(1)})_{(1)}\act a_{(1)})))((g_{(2)}S(h_{(1)})_{(2)})\act a_{(2)})\ot (g_{(3)}S(h_{(1)})_{(3)})h_{(2)}.
\end{align*}
The last equality follows by the  $H$-linearity of the antipode $S$ of $A$. Therefore,
$$
(bS(g_{(1)}\act (S(h_{(1)})_{(1)}\act a_{(1)})))((g_{(2)}S(h_{(1)})_{(2)})\act a_{(2)})\ot (g_{(3)}S(h_{(1)})_{(3)})h_{(2)}=\varepsilon (a)\varepsilon (h)b\ot g.
$$
Next, apply $\id \ot \varepsilon$ to both sides of this equation, set $b=g_{(1)}\act (S(h)_{(1)}\act a_{(1)})$ and replace $a$, $g$, $h$ by $a\sw 2$, $g\sw 2$, $h\sw 2$, respectively, to obtain
\[
((g_{(1)}\act (S(h)_{(1)}\act a_{(1)}))S(g_{(2)}\act (S(h)_{(2)}\act a_{(2)})))((g_{(3)}S(h)_{(3)})\act a_{(3)})
=g\act (S(h)\act a).
\]
Finally, the fact that $A$ is an  $H$-quasimodule coalgebra \eqref{qmcoalg} combined with the antipode property yield
$
g\act (S(h)\act a)=(gS(h))\act a,
$
as required.

 $(2)\Rightarrow (1)$ Assume that condition \eqref{modass} holds. We need to check that $A\lsmash  H$ is a Hopf quasigroup. 

Obviously, $A\lsmash  H$ is a coassociative and counital coalgebra, the unitality of multiplication is a consequence of the unitality of the quasiaction, product and coproduct in Hopf quasigroups. The comultiplication in $A\lsmash  H$  is an algebra homomorphism by the quasimodule coalgebra property \eqref{qmcoalg}  and by \eqref{cocommu}. All this is not different from the case of a smash product of standard Hopf algebras.

It remains to check the Hopf quasigroup identities \eqref{quasi1} and \eqref{quasi2}. They are checked by direct calculations, which we display presently indicating carefully what assumptions and equalities are used at each stage. For all $a,b\in A$ and $g,h\in H$, 
\begin{align*}
S((a\ot h)_{(1)})& ((a\ot h)_{(2)}(b\ot g))
=S(a_{(1)}\ot h_{(1)})((a_{(2)}\ot h_{(2)})(b\ot g))\\
&\!\!\!\!\stackrel{\eqref{smashanti}\eqref{smashmulti}}{=}(S(h_{(1)})_{(1)}\act S(a_{(1)})\ot S(h_{(1)})_{(2)})(a_{(2)}(h_{(2)}\act b)\ot h_{(3)}g)\\
&\stackrel{\eqref{smashmulti}}{=}(S(h_{(1)})_{(1)}\act S(a_{(1)}))(S(h_{(1)})_{(2)}\act (a_{(2)}(h_{(2)}\act b)))\ot S(h_{(1)})_{(3)}(h_{(3)}g)\\
&\stackrel{\eqref{qmalg}}{=}S(h_{(1)})_{(1)}\act (S(a_{(1)})(a_{(2)}(h_{(2)}\act b)))\ot S(h_{(1)})_{(2)}(h_{(3)}g)\\
&\stackrel{\eqref{quasi1}}{=}\varepsilon (a)S(h_{(2)})\act (h_{(3)}\act b)\ot S(h_{(1)})(h_{(4)}g)
\stackrel{\eqref{qm2}\eqref{quasi1}}{=}\varepsilon (a)\varepsilon (h)b\ot g.
\end{align*}
This proves the first of equations \eqref{quasi1}. Next
\begin{align*}
(a&\ot h)_{(1)}(S((a\ot h)_{(2)})(b\ot g))
=(a_{(1)}\ot h_{(1)})(S(a_{(2)}\ot h_{(2)})(b\ot g))\\
&\stackrel{\eqref{smashanti}}{=}(a_{(1)}\ot h_{(1)})(S(h_{(3)})\act S(a_{(2)})\ot S(h_{(2)}))(b\ot g))\\
&\stackrel{\eqref{smashmulti}}{=}(a_{(1)}\ot h_{(1)})(S(h_{(4)})\act S(a_{(2)}))(S(h_{(3)})\act b)\ot S(h_{(2)})g)\\
&\stackrel{\eqref{smashmulti}}{=}(a_{(1)}(h_{(1)}\act ((S(h_{(5)})\act S(a_{(2)}))(S(h_{(4)})\act b))\ot h_{(2)}(S(h_{(3)})g)\\
&\stackrel{\eqref{quasi1}}{=}a_{(1)}(h_{(1)}\act ((S(h_{(3)})\act S(a_{(2)}))(S(h_{(2)})\act b))\ot g
\stackrel{\eqref{qmalg}}{=}a_{(1)}(h_{(1)}\act (S(h_{(2)})\act  (S(a_{(2)})b))\ot g\\
&\stackrel{\eqref{qm2}}{=}\varepsilon (h)a_{(1)}(S(a_{(2)})b)\ot g
\stackrel{\eqref{quasi1}}{=}\varepsilon (a)\varepsilon (h)b\ot g,
\end{align*}
thus proving the second of relations \eqref{quasi1}. It is the proof of  \eqref{quasi2} where the associative law \eqref{modass} is used. The first of identities  \eqref{quasi2}  is proven by the following calculation
\begin{align*}
((b&\ot g)(a\ot h)_{(1)})S((a\ot h)_{(2)})
=((b\ot g)(a_{(1)}\ot h_{(1)}))S(a_{(2)}\ot h_{(2)})\\
&\stackrel{\eqref{smashmulti}\eqref{smashanti}}{=}(b(g_{(1)}\act a_{(1)})\ot g_{(2)}h_{(1)})(S(h_{(3)})\act S(a_{(1)})\ot S(h_{(2)}))\\
&\stackrel{\eqref{smashmulti}}{=}(b(g_{(1)}\act a_{(1)}))((g_{(2)}h_{(1)})\act (S(h_{(4)})\act S(a_{(2)}))\ot (g_{(3)}h_{(2)})S(h_{(3)})\\
&\stackrel{\eqref{quasi2}}{=}(b(g_{(1)}\act a_{(1)}))((g_{(2)}h_{(1)})\act (S(h_{(2)})\act S(a_{(2)}))\ot g_{(3)}\\
&\stackrel{\eqref{modass}}{=}(b(g_{(1)}\act a_{(1)}))(((g_{(2)}h_{(1)})(S(h_{(2)}))\act S(a_{(2)}))\ot g_{(3)}\\
&\stackrel{\eqref{quasi2}}{=}\varepsilon (h)(b(g_{(1)}\act a_{(1)}))(g_{(2)}\act S(a_{(2)}))\ot g_{(3)}
{=}\varepsilon (h)(b(g_{(1)}\act a_{(1)}))S(g_{(2)}\act a_{(2)})\ot g_{(3)}\\
&\stackrel{\eqref{qmcoalg}}{=}\varepsilon (h)(b(g_{(1)}\act a)_{(1)})S(g_{(1)}\act a)_{(2)}\ot g_{(2)}
\stackrel{\eqref{quasi2}}{=}\varepsilon (a)\varepsilon (h)b\ot g,
\end{align*}
where the seventh equality is a consequence of Lemma~\ref{def.qmS}. Finally,
\begin{align*}
((b\ot g)&S((a\ot h)_{(1)}))(a\ot h)_{(2)}
=((b\ot g)S(a_{(1)}\ot h_{(1)}))(a_{(2)}\ot h_{(2)})\\
&\stackrel{\eqref{smashanti}}{=}((b\ot g)(S(h_{(1)})_{(1)}\act S(a_{(1)})\ot S(h_{(1)})_{(2)})(a_{(2)}\ot h_{(2)})\\
&\stackrel{\eqref{smashmulti}}{=} (b(g_{(1)}\act (S(h_{(1)})_{(1)}\act S(a_{(1)})))\ot g_{(2)}S(h_{(1)})_{(2)})(a_{(2)}\ot h_{(2)})\\
&\stackrel{\eqref{smashmulti}}{=} (b(g_{(1)}\act (S(h_{(1)})_{(1)}\act S(a_{(1)})))((g_{(2)}S(h_{(1)})_{(2)})\act a_{(2)})\ot (g_{(3)}S(h_{(1)})_{(3)})h_{(2)}\\
&\stackrel{\eqref{modass}}{=} (b(g_{(1)}S(h_{(1)})_{(1)}\act S(a_{(1)}))((g_{(2)}S(h_{(1)})_{(2)})\act a_{(2)})\ot (g_{(3)}S(h_{(1)})_{(3)})h_{(2)}\\
&\;\;{=}\;\;(b((g_{(1)}S(h_{(1)})_{(1)})_{(1)}\act S(a_{(1)}))((g_{(1)}S(h_{(1)})_{(1)})_{(2)}\act a_{(2)})\ot (g_{(2)}S(h_{(1)})_{(2)})h_{(2)}\\
&\stackrel{\eqref{qmcoalg}}{=}(bS((g_{(1)}S(h_{(1)})_{(1)}\act a)_{(1)})(g_{(1)}S(h_{(1)})_{(1)}\act a)_{(2)}\ot (g_{(2)}S(h_{(1)})_{(2)})h_{(2)}\\
&\stackrel{\eqref{quasi2}}{=}\varepsilon (a)b\ot (gS(h_{(1)}))h_{(2)}=
\varepsilon (a)\varepsilon (h)b\ot g.
\end{align*}
The sixth equality is a consequence of the multiplicativity of the coproduct. Furthermore, Lemma~\ref{def.qmS} is also used in the derivation of the seventh equality. This completes the proof that  $A\lsmash  H$ is a Hopf quasigroup as required.
\endproof

\subsection{Comodule coalgebras  and smash coproducts for Hopf co\-quasi\-groups}

The construction of smash coproducts of Hopf co\-quasi\-groups is dual to the construction described in Section~\ref{subsec.smash}. Thus, first one needs to formulate  the definitions of a quasicomodule algebra and a quasicomodule coalgebra.

\begin{definition}\label{def.coqmalg}
Let $H$ be a Hopf coquasigroup.  An associative and unital algebra $A$ is called a {\em right $H$-quasicomodule algebra} if $(A,\rho^A)$ is a right $H$-quasicomodule and, for all $a, b\in A$, 
\begin{equation}\label{coqmalg}
\rho^A (ab)=\rho^A (a)\rho^A (b), \qquad \rho^A (1)=1\ot 1.
\end{equation}
A counital (but not necessarily coassociative) coalgebra $C$  is termed a {\em right $H$-quasicomodule coalgebra} if $(C,\rho^C)$, $\rho^C: c\mapsto c\suc 0\ot c\suc 1$, is a right $H$-quasicomodule and, for all $c\in C$,
\begin{equation}\label{coqmcoalg}
{c^{(0)}}_{(1)}\ot {c^{(0)}}_{(2)}\ot c^{(1)}= {c_{(1)}}^{(0)}\ot {c_{(2)}}^{(0)}\ot  {c_{(1)}}^{(1)} {c_{(2)}}^{(1)}, \qquad \varepsilon (c^{(0)})c^{(1)}=\varepsilon (c).
\end{equation}
A {\em  right $H$-quasicomodule Hopf coquasigroup} is a Hopf coquasigroup $A$ that is a right $H$-comodule algebra and $H$-comodule coalgebra (by the same coaction).
\end{definition}
Dualising Lemma~\ref{def.qmS} one obtains (see \cite[Lemma~5.13]{KliMaj:Hop})
\begin{lemma}\label{def.coqmS}
The antipode of an $H$-quasicomodule Hopf coquasigroup is $H$-colinear.
\end{lemma}

The necessary conditions for and the construction of a smash coproduct of Hopf coquasigroups are contained in the following (see \cite[Proposition~5.14]{KliMaj:Hop})

\begin{theorem}\label{def.cosmash}
Let $H$ be a Hopf coquasigroup,  $C$ a right $H$-quasicomodule Hopf coquasigroup such that, for all $c\in C$, $h\in H$, 
\begin{equation}\label{commu}
c^{(0)}\ot c^{(1)}h=c^{(0)}\ot hc^{(1)}.
\end{equation}
Then the following statements are equivalent: 
\begin{zlist}
\item
There is a smash coproduct Hopf coquasigroup $C\rcosmash  H$ built on $C\ot H$ with tensor product algebra structure and counit, and with the coproduct and antipode defined by
\begin{equation}\label{smashcomulti}
\Delta (h\ot c)=h_{(1)}\ot {c_{(1)}}^{(0)}\ot h_{(2)}{c_{(1)}}^{(1)}\ot c_{(2)},
\end{equation}
\begin{equation}\label{smashcoanti}
S(h\ot c)=S(hc^{(1)})\ot S(c^{(0)}),
\end{equation}
for all $c\in C, h\in H$.
\item 
For all $c\in A$, 
\begin{equation}\label{comodcoass}
c^{(0)(0)}\ot S(c^{(0)(1)})\ot c^{(1)}=c^{(0)}\ot S({c^{(1)}}_{(1)})\ot {c^{(1)}}_{(2)}.
\end{equation}
\end{zlist}
\end{theorem}

\proof This is dual to Theorem~\ref{def.cosmash}, so we only indicate the key steps leading to and using the coassociative law \eqref{comodcoass}.

$(1)\Rightarrow (2)$ If $C\rcosmash  H$ is a Hopf coquasigroup with the described structure, then necessarily, for all $c\in C$ and $h\in H$, 
\begin{equation}\label{anti}
(h\ot c)_{(1)}S((h\ot c)_{(2)(1)})\ot (h\ot c)_{(2)(2)}
=1\ot 1\ot h\ot c\, ;
\end{equation}
see the first of equations \eqref{coq1}. Develop the left hand side of \eqref{anti} 
\begin{align*}
 (h&\ot c)_{(1)}S((h\ot c)_{(2)(1)})\ot (h\ot c)_{(2)(2)}\\
&\stackrel{\eqref{smashcomulti}}{=} (h_{(1)}\ot {c_{(1)}}^{(0)})S(h_{(2)(1)}{{c_{(1)}}^{(1)}}_{(1)}\ot {c_{(2)(1)}}^{(0)})\ot (h_{(2)(2)}
{{c_{(1)}}^{(1)}}_{(2)}{c_{(2)(1)}}^{(1)}\ot c_{(2)(2)})\\
&\stackrel{\eqref{smashcoanti}}{=} (h_{(1)}\ot {c_{(1)}}^{(0)})(S(h_{(2)(1)}{{c_{(1)}}^{(1)}}_{(1)}{c_{(2)(1)}}^{(0)(1)})\ot S({c_{(2)(1)}}^{(0)(0)})) \\
& \hspace{3in} \ot  (h_{(2)(2)}{{c_{(1)}}^{(1)}}_{(2)}{c_{(2)(1)}}^{(1)}\ot c_{(2)(2)})\\
&\;\;{=}\;\; (h_{(1)}S({{c_{(1)}}^{(1)}}_{(1)}{c_{(2)(1)}}^{(0)(1)})S(h_{(2)(1)})\ot {c_{(1)}}^{(0)}S({c_{(2)(1)}}^{(0)(0)}) \\
& \hspace{3in} \ot  (h_{(2)(2)}{{c_{(1)}}^{(1)}}_{(2)}{c_{(2)(1)}}^{(1)}\ot c_{(2)(2)}),
\end{align*}
where the antimultiplicativity of the antipode is used for the last equality. 
Set $h=1$ in the equality resulting from \eqref{anti} to obtain
\[
S({{c_{(1)}}^{(1)}}_{(1)}{c_{(2)(1)}}^{(0)(1)})\ot {c_{(1)}}^{(0)}S({c_{(2)(1)}}^{(0)(0)})\ot {{c_{(1)}}^{(1)}}_{(2)}{c_{(2)(1)}}^{(1)}\ot c_{(2)(2)}=1\ot 1\ot 1\ot c.
\]
By Lemma~\ref{def.coqmS},
\begin{align*}
S(S({c_{(2)(1)}})^{(0)(1)})S({{c_{(1)}}^{(1)}}_{(1)})\ot {c_{(1)}}^{(0)}S(c_{(2)(1)})^{(0)(0)}\ot {{c_{(1)}}^{(1)}}_{(2)}S({c_{(2)(1)}})^{(1)}&\ot c_{(2)(2)}\\
&=1\ot 1\ot 1\ot c.
\end{align*}
Apply $(1\ot 1\ot 1\ot 1\ot S\ot 1\ot 1)\circ (1\ot 1\ot 1\ot \rho^C\ot 1\ot 1)\circ (1\ot 1\ot 1\ot \rho^C\ot 1)\circ   (1\ot 1\ot 1\ot \Delta)$ to the equation above to conclude that
\begin{align*}
S&(S({c_{(2)(1)}})^{(0)(1)})S({{c_{(1)}}^{(1)}}_{(1)})\ot {c_{(1)}}^{(0)}S(c_{(2)(1)})^{(0)(0)}\ot {{c_{(1)}}^{(1)}}_{(2)}S({c_{(2)(1)}})^{(1)}\ot 
{c_{(2)(2)(1)}}^{(0)(0)}\\
&\ot S({c_{(2)(2)(1)}}^{(0)(1)})\ot {c_{(2)(2)(1)}}^{(1)}\ot c_{(2)(2)(2)}
=1\ot 1\ot 1\ot {c_{(1)}}^{(0)(0)}\ot S({c_{(1)}}^{(0)(1)})\ot {c_{(1)}}^{(1)}\ot c_{(2)}.
\end{align*}
Next, multiply the first component by the fifth one, the second component  by the fourth one and the third component by the sixth one to obtain
\begin{align*}
S&({c_{(2)(2)(1)}}^{(0)(1)})S(S({c_{(2)(1)}})^{(0)(1)})S({{c_{(1)}}^{(1)}}_{(1)})\ot {c_{(1)}}^{(0)}S(c_{(2)(1)})^{(0)(0)}{c_{(2)(2)(1)}}^{(0)(0)} \\
&\ot {{c_{(1)}}^{(1)}}_{(2)}S({c_{(2)(1)}})^{(1)}{c_{(2)(2)(1)}}^{(1)}\ot c_{(2)(2)(2)}=S({c_{(1)}}^{(0)(1)})\ot {c_{(1)}}^{(0)(0)}\ot {c_{(1)}}^{(1)}\ot c_{(2)}.
\end{align*}
The next step requires the use of the antipode property and the quasicomodule algebra condition \eqref{coqmalg}
\begin{align*}
S(S&({c_{(2)(1)}}){c_{(2)(2)(1)}})^{(0)(1)}S({{c_{(1)}}^{(1)}}_{(1)})\ot {c_{(1)}}^{(0)}(S(c_{(2)(1)}){c_{(2)(2)(1)}})^{(0)(0)} \\
& \ot {{c_{(1)}}^{(1)}}_{(2)}(S({c_{(2)(1)}}){c_{(2)(2)(1)}})^{(1)}\ot c_{(2)(2)(2)}=S({c_{(1)}}^{(0)(1)})\ot {c_{(1)}}^{(0)(0)}\ot {c_{(1)}}^{(1)}\ot c_{(2)}.
\end{align*}
The coassociative law \eqref{comodcoass} follows by the use of coquasigroup identity \eqref{coq1} and  application of $\varepsilon$ to the fourth component.
This completes the proof of the necessity of \eqref{comodcoass}.

$(2)\Rightarrow (1)$ We only display the verification of  coquasigroup identities \eqref{coq1} as this is where \eqref{comodcoass} is used. First,
\begin{align*}
S&((h\ot c)_{(1)})(h\ot c)_{(2)(1)}\ot (h\ot c)_{(2)(2)}\\
&\stackrel{\eqref{smashcomulti}}{=}S(h_{(1)}\ot {c_{(1)}}^{(0)})(h_{(2)(1)}{{c_{(1)}}^{(1)}}_{(1)}\ot {c_{(2)(1)}}^{(0)})\ot h_{(2)(2)}
{{c_{(1)}}^{(1)}}_{(2)}{c_{(2)(1)}}^{(1)}\ot c_{(2)(2)}\\
&\stackrel{\eqref{smashcoanti}}{=} (S(h_{(1)}{c_{(1)}}^{(0)(1)})\ot S({c_{(1)}}^{(0)(0)}))(h_{(2)(1)}{{c_{(1)}}^{(1)}}_{(1)}\ot {c_{(2)(1)}}^{(0)}) \\
&\hspace{3in} \ot h_{(2)(2)}{{c_{(1)}}^{(1)}}_{(2)}{c_{(2)(1)}}^{(1)}\ot c_{(2)(2)}\\
&\;\;=\;\; S({c_{(1)}}^{(0)(1)})S(h_{(1)})h_{(2)(1)}{{c_{(1)}}^{(1)}}_{(1)}\ot S({c_{(1)}}^{(0)(0)}) {c_{(2)(1)}}^{(0)} \\
& \hspace{3in} \ot h_{(2)(2)}{{c_{(1)}}^{(1)}}_{(2)}{c_{(2)(1)}}^{(1)}\ot c_{(2)(2)}\\
&\stackrel{\eqref{coq1}}{=} S({c_{(1)}}^{(0)(1)}){{c_{(1)}}^{(1)}}_{(1)}\ot S({c_{(1)}}^{(0)(0)}) {c_{(2)(1)}}^{(0)} \ot h{{c_{(1)}}^{(1)}}_{(2)}{c_{(2)(1)}}^{(1)}\ot c_{(2)(2)}\\
&\stackrel{\eqref{comodcoass}}{=} \! S({{c_{(1)}}^{(1)}}_{(1)}){{c_{(1)}}^{(1)}}_{(2)(1)}\ot S({c_{(1)}}^{(0)}){c_{(2)(1)}}^{(0)}\ot h{{c_{(1)}}^{(1)}}_{(2)(1)}{c_{(2)(2)}}^{(1)}{c_{(2)(1)}}^{(1)}\ot c_{(2)(2)}\\
&\stackrel{\eqref{coq1}}{=} 1\ot S({c_{(1)}}^{(0)}){c_{(2)(1)}}^{(0)}\ot h{{c_{(1)}}^{(1)}}{c_{(2)(2)}}^{(1)}{c_{(2)(1)}}^{(1)}\ot c_{(2)(2)}\\
&\;\,=\;\, 1\ot S({c_{(1)}})^{(0)}{c_{(2)(1)}}^{(0)}\ot hS({{c_{(1)}})^{(1)}}{c_{(2)(2)}}^{(1)}{c_{(2)(1)}}^{(1)}\ot c_{(2)(2)}\\
&\stackrel{\eqref{coqmalg}}{=} 1\ot (S(c_{(1)})c_{(2)(1)})^{(0)}\ot h(S(c_{(1)})c_{(2)(2)})^{(1)}{c_{(2)(1)}}^{(1)}\ot c_{(2)(2)}
\stackrel{\eqref{coq1}}{=}1\ot 1\ot h\ot c,
\end{align*}
where the third equality is a consequence of the antimultiplicativity of the antipode, and the seventh equality is a consequence of Lemma~\ref{def.coqmS}. Second,
\begin{align*}
  (h&\ot c)_{(1)}S((h\ot c)_{(2)(1)})\ot (h\ot c)_{(2)(2)}\\
&\stackrel{\eqref{smashcomulti}}{=} (h_{(1)}\ot {c_{(1)}}^{(0)})S(h_{(2)(1)}{{c_{(1)}}^{(1)}}_{(1)}\ot {c_{(2)(1)}}^{(0)})\ot h_{(2)(2)}
{{c_{(1)}}^{(1)}}_{(2)}{c_{(2)(1)}}^{(1)}\ot c_{(2)(2)}\\
&\stackrel{\eqref{smashcoanti}}{=} (h_{(1)}\ot {c_{(1)}}^{(0)})(S(h_{(2)(1)}{{c_{(1)}}^{(1)}}_{(1)}{c_{(2)(1)}}^{(0)(1)})\ot S({c_{(2)(1)}}^{(0)(0)})) \\
& \hspace{3in} \ot h_{(2)(2)}{{c_{(1)}}^{(1)}}_{(2)}{c_{(2)(1)}}^{(1)}\ot c_{(2)(2)}\\
&\;\;{=}\;\; h_{(1)}S({{c_{(1)}}^{(1)}}_{(1)}{c_{(2)(1)}}^{(0)(1)})S(h_{(2)(1)})\ot {c_{(1)}}^{(0)}S({c_{(2)(1)}}^{(0)(0)}) \\
&\hspace{3in} \ot  h_{(2)(2)}{{c_{(1)}}^{(1)}}_{(2)}{c_{(2)(1)}}^{(1)}\ot c_{(2)(2)}\\
&\stackrel{\eqref{comodcoass}}{=}h_{(1)}S({{c_{(1)}}^{(1)}}_{(1)}{{c_{(2)(1)}}^{(1)}}_{(1)})S(h_{(2)(1)})\ot {c_{(1)}}^{(0)}S({c_{(2)(1)}}^{(0)}) \\
& \hspace{3in} \ot  h_{(2)(2)}{{c_{(1)}}^{(1)}}_{(2)}{{c_{(2)(1)}}^{(1)}}_{(2)}\ot c_{(2)(2)}\\
&\;\;{=}\;\;h_{(1)}S({{c_{(1)}}^{(1)}}_{(1)}{{S(c_{(2)(1)})}^{(1)}}_{(1)})S(h_{(2)(1)})\ot {c_{(1)}}^{(0)}S(c_{(2)(1)})^{(0)} \\
& \hspace{3in} \ot  h_{(2)(2)}{{c_{(1)}}^{(1)}}_{(2)}{{S(c_{(2)(1)})}^{(1)}}_{(2)}\ot c_{(2)(2)}\\
&\stackrel{\eqref{coqmalg}}{=} h_{(1)}S(c_{(1)}{S(c_{(2)(1)}))^{(1)}}_{(1)}S(h_{(2)(1)})\ot (c_{(1)}S(c_{(2)(1)}))^{(0)} \\
&\hspace{3in} \ot h_{(2)(2)}(c_{(1)}{S(c_{(2)(1)}))^{(1)}}_{(2)}\ot c_{(2)(2)}\\
&\stackrel{\eqref{coq1}}{=}h_{(1)}S(h_{(2)(1)})\ot 1\ot h_{(2)(2)}\ot c
\stackrel{\eqref{coq1}}{=}1\ot 1\ot h\ot c,
\end{align*}
where again the third equality is a consequence of the antimultiplicativity of the antipode, and the fifth one is a consequence of Lemma~\ref{def.coqmS}. 
\endproof

 \section*{Acknowledgements} 
 The research of Zhengming Jiao is supported by the Natural Science Foundation of Henan Province (grant no.\ 102300410049).

 \end{document}